\documentclass[a4paper,10pt]{article}

\usepackage{amsmath}
\usepackage{amssymb}
\usepackage{amsthm}
\usepackage{mathrsfs}
\usepackage{color}
\usepackage[dvips]{graphicx}
\usepackage[all]{xy}

\theoremstyle{plain}
\newtheorem{them}{Theorem}[section]
\newtheorem{lemma}[them]{Lemma}
\newtheorem{prop}[them]{Proposition}

\theoremstyle{definition}
\newtheorem{defi}[them]{Definition}
\newtheorem{exam}[them]{Example}

\newtheorem*{mthem}{Main Theorem}

\newcommand{\lmap}[3]{#1:#2 \longrightarrow #3}
\newcommand{\map}[3]{#1:#2 \rightarrow #3}

\begin{document}

\title{Ramified coverings of small categories}
\author{Kazunori Noguchi \thanks{noguchi@math.shinshu-u.ac.jp}}
\date{}
\maketitle
\begin{abstract}
We introduce a ramified covering of small categories, and we show three properties of the notion: the Riemann-Hurwitz formula holds for a ramified covering of finite categories, the zeta function of $C$ divides that of $\widetilde{C}$ for a ramified covering $\map{P}{\widetilde{C}}{C}$ of finite categories, and the classifying space of a $d$-fold ramified covering of small categories is also a $d$-fold ramified covering in the sense of Dold \cite{Dol86}.
\end{abstract}


\thispagestyle{empty}

\section{Introduction}

A covering is an interesting and important tool for geometry. For example, a covering space is used for computations of fundamental groups and it has an analogy of Galois theory, see, for example, \cite{Hat02} and \cite{May99}. A covering space should be called ''unramified covering". A ramified covering for topological spaces is defined by Smith \cite{Smi83} and Dold \cite{Dol86}, but an well-known example of ramified coverings would be the one for Riemann surfaces. A ramified covering has important properties as same as unramified coverings do, for instance, the Riemann-Hurwitz formula holds, where it states a relationship between the Euler characteristics of a total space and a base space.

In this paper, we define a ramified covering of small categories. An unramified covering of small categories has already been defined, and several authors have studied about it. By the works of Bridson and Haefliger, we can find many important properties of unramified coverings in \cite{BH99}, for example, the monodromy theorem, the path lifting theorem and so on. May studied about unramified coverings of groupoids \cite{May99}. Tanaka defined a model structure on the categories of small categories, called $1$-type model structure \cite{Tan}. An unramified covering is a fibration in the sense of $1$-type model structure. Cibils and MacQuarrie studied about Galois coverings of small categories \cite{CM}.

Ramified coverings of small categories in this paper satisfy many desirable properties.

\begin{mthem}
Let $\map{P}{\widetilde{C}}{C}$ be a $d$-fold ramified covering of small categories. Then, we obtain the following results:
\begin{enumerate}
\item Suppose $\widetilde{C}$ and $C$ are finite categories. The category $\widetilde{C}$ has series Euler characteristic if and only if $C$ has series Euler characteristic. In this case,
$$\chi_{\sum}(\widetilde{C})=d\chi_{\sum}(C)-V,$$
where $$V=\sum_{\widetilde{x}\in \mathrm{Ob}(\widetilde{C})} (e(\tilde{x})-1),$$
and $e(\tilde{x})$ is the ramification number of $\tilde{x}$.
\item Suppose $\widetilde{C}$ and $C$ are finite categories. The zeta function of $C$ divides that of $\widetilde{C}$, that is,
$$\zeta_{\widetilde{C}}(z)=\zeta_C(z)^d (1-z)^V.$$
\item The map $\map{BP}{B\widetilde{C}}{BC}$ is a $d$-fold ramified covering in the sense of Dold \cite{Dol86}, where $B$ is the classifying space functor.
\end{enumerate}
\end{mthem}

The first result is an analogue of the Riemann-Hurwitz formula for Riemann surfaces. For a ramified covering $\map{p}{\widetilde{X}}{X}$ of Riemann surfaces under certain condition, the following is well-known as the Riemann-Hurwitz formula: $$\chi(\widetilde{X})=d \chi(X) -V,$$
where $d$ is the degree of $p$ and $$V=\sum_{\tilde{x} \in \widetilde{X}} (e(\tilde{x})-1).$$

Euler characteristic for categories is defined in various ways. Leinster defined the \textit{Euler characteristic of a finite category} in \cite{Leia}. This is the first Euler characteristic for categories, and later several authors defined, the \textit{series Euler characteristic} by Berger-Leinster \cite{Leib}, the \textit{$L^2$-Euler characteristic} by Fiore-L\"uck-Sauer \cite{FLS}, the \textit{extended $L^2$-Euler characteristic} \cite{Nog} and the \textit{Euler characteristic of $\mathbb{N}$-filtered acyclic categories} by the author \cite{Nog11}. See \cite{Nog} for relationships among them. In this paper, we only show the Riemann-Hurwitz formula holds for the series Euler characteristic. The author does not know if the other Euler characteristics have such property or not.

Graph theoretic analogue of Riemann-Hurwitz formula is considered in \cite{BN09}.

The second result is a generalization of Theorem 4.5 of \cite{NogA}. Graph theoretic analogue of this result is also considered in \cite{MM10}, \cite{ST96} and \cite{Ter11}.

This result is a categorical analogue of the Dedekind conjecture that states if $K_1$ and $K_2$ are number fields and $K_1 \subset K_2$, then the Dedekind zeta function of $K_1$ divides that of $K_2$. A covering of small categories is an analogy of Galois theory as same as a covering of topological spaces is so (see, for example, \cite{Hat02} and \cite{May99}). Fundamental theorem of Galois theory is if $K/F$ is a finite Galois extension, the set of intermediate fields of $K$ and $F$ is bijective to the set of subgroups of the Galois group $\mathrm{Gal}(K/F)$ $$\xymatrix{K\ar@{-}[d]\ar@{<->}[r]&\{e\} \ar@{}[d]|-{\bigcap}\\ L\ar@{-}[d] \ar@{<->}[r]^{1 :1}&H \ar@{}[d]|-{\bigcap}\\ F\ar@{<->}[r]&\mathrm{Gal}(K/F).}$$ For a covering of small categories $\map{\tilde{P}}{\tilde{E}}{B}$ where $\widetilde{E}$ is the universal covering of $B$, the set of the isomorphism classes of intermediate coverings of $\tilde{P}$ is bijective to the set of subgroups of the fundamental group $\pi_1(B)$
$$\xymatrix{\widetilde{E}\ar@/_/[dd]_{\tilde{P}}\ar[d]\ar@{<->}[r]&\{e\} \ar@{}[d]|-{\bigcap}\\ E\ar[d] \ar@{<->}[r]^{1 :1}&H \ar@{}[d]|-{\bigcap}\\ F\ar@{<->}[r]&\pi_1(B)}$$ 
(see Corollary 2.24 of \cite{Tan}). We have the following correspondences:
\begin{eqnarray*}
\text{coverings}&\leftrightarrow&\text{extension of fields} \\ \pi_1 & \leftrightarrow&\text{Galois groups} \\ \text{intermediate coverings}& \leftrightarrow&\text{intermediate fields} \\ \vdots &&\vdots
\end{eqnarray*}
For an analogy between coverings of spaces and extensions of fields, see \cite{Mor12}. By the diagrams above, we can conclude that the relationship between zeta functions and coverings is an analogue of the Dedekind conjecture.

When $\map{P}{\widetilde{C}}{C}$ is an unramified covering, it is known that $BP$ is a covering space, see,  for instance, \cite{Tan}. The third result is a generalization of such fact when the fiber of $P$ is finite. Smith defined ramified coverings of spaces \cite{Smi83}, and later Dold simplified and developed Smith's definition and theory \cite{Dol86}. We use Dold's definition in this paper.

\section{Ramified coverings of small categories}

\subsection{Notation and terminology}

Before we introduce a ramified covering of small categories, let us recall unramified coverings of small categories \cite{BH99}.

Let $C$ be a small category. For an object $x$ of $C$, let $S(x)$ be the set of morphisms of $C$ whose source is $x$
$$S(x)=\{\map{f}{x}{*}\in \mathrm{Mor}(C)\},$$
and $T(x)$ is the set of morphisms of $C$ whose target is $x$
$$T(x)=\{\map{g}{*}{x}\in \mathrm{Mor}(C)\}.$$ 
We denote by $\overline{S}(x)$ $S(x)-\{1_x\}$.

A category $C$ is \textit{connected} if $C$ is a non-empty category and there exists a zig-zag sequence of morphisms in $C$
$$\xymatrix{x\ar[r]^{f_1}&x_1&x_2\ar[r]^{f_3}\ar[l]_{f_2}&\dots&y\ar[l]_(0.4){f_n}}$$
for any objects $x$ and $y$ of $C$. We do not have to care about the direction of the last morphism $f_n$ since we can insert an identity morphism to the sequence.

A functor $\map{P}{\widetilde{C}}{C}$ is an \textit{unramified covering} if the following two restrictions of $P$
$$P:S(\tilde{x})\longrightarrow S(x)$$ $$P:T(\tilde{x})\longrightarrow T(x)$$
are bijections for any object $\tilde{x}$ of $\widetilde{C}$ and $P(\tilde{x})=x$. This condition is an analogue of the condition of an unramified covering of graphs (see \cite{ST96}).

Let
$$N_n(C)=\{\xymatrix{(x_0\ar[r]^{f_1}&x_1\ar[r]^{f_2}&\dots\ar[r]^{f_n}&x_n)} \text{ in } C \}$$
and
$$\overline{N_n}(C)=\{\xymatrix{(x_0\ar[r]^{f_1}&x_1\ar[r]^{f_2}&\dots\ar[r]^{f_n}&x_n)} \text{ in } C \mid  f_i \not = 1\}.$$
The difference between them is just one thing that identity morphisms are used or not. For $m=0$, we set $N_0(C)=\overline{N_0}(C)=\mathrm{Ob}(C)$.

\subsection{Definition}

\begin{defi}
Suppose $\map{P}{\widetilde{C}}{C}$ is a functor and $C$ is connected. Then, $P$ is a \textit{ramified covering} if it satisfies the following conditions:
\begin{enumerate}
\item For each object $\tilde{x}$ of $\widetilde{C}$, there exists a natural number $e(\tilde{x})$, called the \textit{ramification number} of $P$ at $\tilde{x}$.
\item The map $\map{P}{T(\tilde{x})}{T(x)}$ is a bijection for any object $\tilde{x}$ of $\widetilde{C}$ and $P(\tilde{x})=x$.
\item The map $\map{P}{\overline{S}(\tilde{x})}{\overline{S}(x)}$ is an $e(\tilde{x})$ to one map for any object $\tilde{x}$ of $\widetilde{C}$ and $P(\tilde{x})=x$.
\end{enumerate}
If $\overline{S}(\tilde{x})$ and  $\overline{S}(x)$ are both empty-sets, we regard $\map{P}{\emptyset}{\emptyset}$ as a one-to-one map, so that the ramification number $e(\tilde{x})$ is $1$. Since $e(\tilde{x})\ge 1$, we do not allow the case that $\overline{S}(\tilde{x})=\emptyset$ and $S(x)\not = \emptyset$.
\end{defi}

\begin{exam}
We introduce only two simple examples here, but the other examples will be given in \S \ref{Exam}.
\begin{enumerate}
\item An unramified covering $\map{P}{\widetilde{C}}{C}$ is a ramified covering as $e(\tilde{x})=1$ for any object $\tilde{x}$ of $\mathrm{Ob}(\widetilde{C})$.
\item Let $$\widetilde{C}=\xymatrix{\tilde{y}_1&&\tilde{y}_2\\ &\tilde{x}\ar[ur]^{\tilde{f}_2} \ar[ul]_{\tilde{f}_1}}$$ and $C=\xymatrix{x\ar[r]^f&y}$. Define a functor $\map{P}{\widetilde{C}}{C}$ by eliminating the tildes and the indexes, for instance, $P(\tilde{y}_1)=y$. Then, $\map{P}{\widetilde{C}}{C}$ is a ramified covering, but not an unramified covering.
\end{enumerate}
\end{exam}

\subsection{Preparation}

In this subsection, we prove lemmas needed later.

\begin{lemma}\label{iff}
Suppose $\map{P}{\widetilde{C}}{C}$ is a ramified covering and $\tilde{f}:\tilde{x}\rightarrow \tilde{y}$ is a morphism in $\widetilde{C}$. Then, $P(\tilde{f})=\map{f}{x}{y}$ is an identity morphism if and only if $\tilde{f}$ is an identity morphism.
\begin{proof}
The bijection $\map{P}{T(\tilde{y})}{T(y)}$ implies this fact.
\end{proof}
\end{lemma}

A functor $\map{F}{C}{D}$ is \textit{finite to one} if for any object $x$ of $D$ and any morphism $f$ of $D$, their inverse images by $P$ are finite sets, that is,
$$\# P^{-1}(x)=\{y\in \mathrm{Ob}(C) \mid P(y)=x\} < +\infty, $$
$$\# P^{-1}(f)=\{g\in \mathrm{Mor}(C) \mid P(g)=f\} <+ \infty. $$

\begin{lemma}\label{cycle}
Suppose $\map{P}{\widetilde{C}}{C}$ is finite to one and a ramified covering, and there exists the following sequence:
$$\xymatrix{\widetilde{\mathbf{f}} =\tilde{x}_0\ar[r]^(0.6){\tilde{f}_1} &\tilde{x}_1\ar[r]^{\tilde{f}_2}&\cdots\ar[r]^{\tilde{f}_{n-1}} &\tilde{x}_{n-1}\ar[r]^(0.6){\tilde{f}_n} &\tilde{x}_0 }$$
such that each $\widetilde{f_i}$ is not an identity morphism. Then, $e(\tilde{x}_0)=1$. Hence, all the $e(\tilde{x}_i)$ are $1$. 
\begin{proof}
Suppose $e(\tilde{x}_0)\ge 2$. Let $$P(\widetilde{\mathbf{f}})=\xymatrix{x_0\ar[r]^{f_1}&x_1\ar[r]^{f_2}&\cdots\ar[r]^{f_{n-1}}&x_{n-1}\ar[r]^{f_{n}}&x_0}.$$Then, there exists a morphism $\map{\tilde{f}_{1,1}}{\tilde{x}_0}{\tilde{x}_{1,1}}$ in $\widetilde{C}$ such that $P(\tilde{f}_{1,1})=f_1$ and $\tilde{f}_{1,1}\not = \tilde{f}_1$. If $\tilde{x}_1=\tilde{x}_{1,1}$, then this fact contradicts to the bijectivity of the map $\map{P}{T(\tilde{x}_1)}{T(x_1)}$, so that $\tilde{x}_1\not =\tilde{x}_{1,1}$. Since $e(\tilde{x}_i)\ge 1$, we have the following diagram
$$\xymatrix{\tilde{\mathbf{f}} =\tilde{x}_0\ar[dr]^{\tilde{f}_{1,1}}\ar[r]^{\tilde{f}_1} &\tilde{x}_1\ar[r]^{\widetilde{f_2}}&\tilde{x}_2\ar[r]&\cdots\ar[r]^{\tilde{f}_{n-1}} &\tilde{x}_{n-1}\ar[r]^{\tilde{f}_n} &\tilde{x}_0\\
&\tilde{x}_{1,1}\ar[r]^{\tilde{f}_{1,2}}&\tilde{x}_{1,2}\ar[r]^{\tilde{f}_{1,3}}&\cdots\ar[r]^{\tilde{f}_{1,n-1}}&\tilde{x}_{1,n-1}\ar[r]^{\tilde{f}_{1,n}}&\tilde{x}_{1,0}\ar[r]^{\tilde{f}_{2,1}}&\tilde{x}_{2,1}
 }$$
such that for each $i$, $\tilde{x}_i\not = \tilde{x}_{1,i}$ and $P(\tilde{f_i})=P(\tilde{f}_{1,i})=f_i$. Moreover, $\tilde{x}_1,\tilde{x}_{1,1}$ and $\tilde{x}_{2,1}$ are all distinct, since if two of them are the same (put it $\tilde{z}$ now), this contradicts to the bijectivity of the map $\map{P}{T(\tilde{z})}{T(x_1)}$. By repeating this process, we obtain infinitely many objects $\tilde{x},\tilde{x}_{1,1},\tilde{x}_{2,1},\dots$ that lie above $x_1$, and this contradicts that $P$ is finite to one. Hence, $e(\tilde{x}_0)=1$.
\end{proof}
\end{lemma}

\begin{lemma}\label{y-lemma}
Suppose $\map{P}{\widetilde{C}}{C}$ is a ramified covering and $\map{\widetilde{f}}{\widetilde{x}}{\tilde{y}}$ is non-identity morphism in $\widetilde{C}$. Then, the ramification number of $\tilde{y}$ is $1$.
\begin{proof}
Suppose $e(\tilde{y})\ge 2$. Then, $\overline{S}(\tilde{y})$ and $\overline{S}(y)$ are not empty, since if the both are empty, then $e(\tilde{y})=1$. Let $P(\tilde{f})=\map{f}{x}{y}$. We can take a morphism $\map{g}{y}{z}$ of $\overline{S}(y)$, and there exist at least two morphisms $\map{\tilde{g}_1}{\tilde{y}}{\tilde{z}_1}$ and $\map{\tilde{g}_2}{\tilde{y}}{\tilde{z}_2}$ such that $\tilde{g}_1 \not =\tilde{g}_2$ and $P(\tilde{g}_1) =P(\tilde{g}_2)=g$. Then, we can show $\tilde{z}_1\not =\tilde{z}_2$ as in the previous proof. Hence, the morphisms $\tilde{g}_1\circ \tilde{f}$ and $\tilde{g}_2\circ \tilde{f}$ are distinct, but $P(\tilde{g}_1\circ \tilde{f})$ and $P(\tilde{g}_2\circ \tilde{f})$ are the same morphism. Then, $e(\tilde{x})=1$, since if $e(\tilde{x})\ge 2$, then at least $e(\tilde{x})+1$ morphisms lie above $g\circ f$. Since $e(\tilde{x})=1$, one of them is an identity morphism. If $\tilde{g}_1\circ \tilde{f}$ is an identity morphism, then $\tilde{x}=\tilde{z}_1$. Lemma \ref{cycle} implies $e(\tilde{y})=1$, so this contradicts to the assumption given at the beginning of this proof. Hence, we obtain $e(\tilde{y})=1$.
\end{proof}
\end{lemma}

Suppose $C$ is a small category and $x$ is an object of $C$. Then, let
$$N_n(C)^x=\{ \xymatrix{x_0\ar[r]^{f_1}&x_1\ar[r]^{f_2}&\cdots\ar[r]^{f_{n}}&x_{n}} \in N_n(C) \mid x_0=x\}.$$
We also define $\overline{N_n}(C)^x$ in similar way.

\begin{lemma}\label{non-degenerate}
Suppose $\map{P}{\widetilde{C}}{C}$ is a ramified covering and $\tilde{x}$ is an object of $\widetilde{C}$ and $P(\tilde{x})=x$. Then, the map
$$\map{P}{\overline{N_n}(\widetilde{C})^{\tilde{x}}}{\overline{N_n}(C)^x}$$
is an $e(\tilde{x})$ to one map for any $n\ge 1$.
\begin{proof}
Given $$\mathbf{f}=\xymatrix{x\ar[r]^{f_1}&x_1\ar[r]^{f_2}&\cdots\ar[r]^{f_{n}}&x_{n}}$$
of $\overline{N_n}(C)^x$, there exist $e(\tilde{x})$ morphisms $\map{\tilde{f}_{i,1}}{\tilde{x}}{\tilde{x}_{i,1}}$ such that $P(\tilde{f}_{i,1})=f_1$. Lemma \ref{y-lemma} implies that there exactly exist $e(\tilde{x})$ lifts of $\mathbf{f}$. Hence, $P$ is $e(\tilde{x})$ to one.
\end{proof}
\end{lemma}

\begin{lemma}\label{degenerate}
Suppose $\map{P}{\widetilde{C}}{C}$ is a ramified covering and $\tilde{x}$ is an object of $\widetilde{C}$ and $P(\tilde{x})=x$. Then, the map
$$\map{P}{N_n(\widetilde{C})^{\tilde{x}}}{N_n(C)^x}$$
is an $e(\tilde{x})$ to one map except for $\mathbf{1}_x$ for any $n\ge 1$, where $\mathbf{1}_{\tilde{x}}$ consists of only the identity morphism $1_{\tilde{x}}$.
\begin{proof}
This is proved in similar way of the proposition above. Note that only the element $\mathbf{1}_{\tilde{x}}$ has one-to-one correspondence to $\mathbf{1}_x$ by Lemma \ref{iff}.
\end{proof}

\end{lemma}

\begin{prop}\label{degree}
Let $\map{P}{\widetilde{C}}{C}$ be a ramified covering. For an object $x$ of $C$, let $R(x)$ be the set of objects of $\widetilde{C}$ that belong to $P^{-1}(x)$, where each $\tilde{x}$ of $P^{-1}(x)$ occurs $e(\tilde{x})$ times. Then, the cardinality of $R(x)$ does not depend on the choice of $x$.
\begin{proof}
It suffices to show that if there exists a morphism $\map{f}{x}{y}$ in $C$, then $R(x)\cong R(y)$. If it is shown, for any object $x$ and $y$ of $C$, there exists the following sequence:
$$\xymatrix{x\ar[r]^{f_1}&x_1&x_2\ar[r]^{f_3}\ar[l]_{f_2}&\dots&y\ar[l]_(0.4){f_n}},$$
then we obtain $$R(x)\cong R(x_1)\cong \cdots \cong R(y).$$

Let $\map{f}{x}{y}$ be a morphism in $C$. Then, for an object $\tilde{x}$ of $R(x)$, there exist $e(\tilde{x})$ morphisms $\map{\tilde{f}_i}{\tilde{x}}{\tilde{y}_i}$ for $1\le i\le e(\tilde{x})$ such that $P(\tilde{f}_i)=f$. We label $\map{\tilde{f}_i}{\tilde{x}}{\tilde{y}_i}$ by $\map{\tilde{f}_i}{\tilde{x}_i}{\tilde{y}_i}$. All of $\tilde{x}_i$ are, in fact, the same object, but they are distinguished in $R(x)$. Define a map $\map{\varphi_f}{R(x)}{R(y)}$ by $\varphi_f(\tilde{x}_i)=\tilde{y}_i$. It follows from the bijectivity with respect to target that $\varphi_f$ is bijective.
\end{proof}
\end{prop}

\begin{defi}
Let $\map{P}{\widetilde{C}}{C}$ be a ramified covering. Define the \textit{degree} of $P$ by the cardinality of $R(x)$ for an object $x$ of $C$. Proposition \ref{degree} implies this definition is well-defined. 

For a natural number $d$, $P$ is a \textit{$d$-fold ramified covering} if $\# R(x)=d$.
\end{defi}

\subsection{Proof of main theorem}

In this subsection, we give a proof of our main theorem.

 A finite category $C$ has \textit{series Euler characteristic} if we can substitute $t=-1$ in the rational function

$$ \frac{\text{sum}(\text{adj}(I-(A_C-I)t ))}{|I-(A_C-I)t |}, $$ 
where $I$ is the unit matrix and $A_C$ is the adjacency matrix of $C$ \cite{Leib}. In this case, the \textit{series Euler characteristic} of $C$ is defined by the value of the rational function at $t=-1$. The rational function is the analytic continuation of the power series $\sum^{\infty}_{n=0} \# \overline{N_n}(C) t^n$ (Theorem 2.2 of \cite{Leib}).

\begin{them}[Riemann-Hurwitz for categories]
Suppose $\map{P}{\widetilde{C}}{C}$ is a $d$-fold ramified covering of finite categories. Then, $\widetilde{C}$ has series Euler characteristic if and only if $C$ has series Euler characteristic. In this case,
$$\chi_{\sum}(\widetilde{C})=d\chi_{\sum}(C)-V,$$
where $$V=\sum_{\widetilde{x}\in \mathrm{Ob}(\widetilde{C})} (e(\tilde{x})-1).$$
\begin{proof}
Let $\chi_C(t)=\sum^{\infty}_{n=0}\# \overline{N_n}(C)t^n$. Then, Lemma \ref{non-degenerate} implies
\begin{eqnarray*}
\chi_{\widetilde{C}}(t)&=&\sum^{\infty}_{n=0}\# \overline{N_n}(\widetilde{C})t^n \\
&=&\# N_0(\widetilde{C})+\sum^{\infty}_{n=1}\# \overline{N_n}(\widetilde{C})t^n  \\
&=&\# N_0(\widetilde{C})+\sum_{x\in \mathrm{Ob}(C)}\sum_{\tilde{x}\in P^{-1}(x)} \sum^{\infty}_{n=1}\# \overline{N_n}(\widetilde{C})^{\tilde{x}} t^n \\
&=&\# N_0(\widetilde{C})+\sum_{x\in \mathrm{Ob}(C)} \bigg(\bigg(\sum_{\tilde{x}\in P^{-1}(x)}e(\tilde{x})\bigg) \sum^{\infty}_{n=1}\# \overline{N_n}(C)^x t^n\bigg) \\
&=&d\chi_C(t)-d\# N_0(C)+\# N_0(\widetilde{C}).
\end{eqnarray*}

So $\widetilde{C}$ has series Euler characteristic if and only if we can substitute $t=-1$ in $$d \frac{\text{sum}(\text{adj}(I-(A_C-I)t ))}{|I-(A_C-I)t |} -d\# N_0(C)+\# N_0(\widetilde{C})$$ if and only if we can substitute $t=-1$ in $$ \frac{\text{sum}(\text{adj}(I-(A_C-I)t ))}{|I-(A_C-I)t |}$$ if and only if $C$ has series Euler characteristic. Hence, the first claim is proven. 

 If $\widetilde{C}$ has series Euler characteristic, then we have
$$\chi_{\sum}(\widetilde{C})=d\chi_{\sum}(C)-d\# N_0(C)+\# N_0(\widetilde{C}).$$
Here we have
\begin{eqnarray*}
V&=&\sum_{\tilde{x}\in \mathrm{Ob}(\widetilde{C})}(e(\tilde{x})-1) \\
&=&\sum_{\tilde{x}\in \mathrm{Ob}(\widetilde{C})}e(\tilde{x})-\# N_0(\widetilde{C}) \\
&=&\sum_{x\in \mathrm{Ob}(C)} \sum_{\tilde{x}\in P^{-1}(x)}e(\tilde{x})-\# N_0(\widetilde{C}) \\
&=&d \# N_0(C)-\# N_0(\widetilde{C}).
\end{eqnarray*}
Hence, we obtain the result.
\end{proof}
\end{them}

Let $C$ be a finite category. Then, the \textit{zeta function} $\zeta_C(z)$ of $C$ is defined by
$$\zeta_C(z)=\exp\left( \sum_{m=1}^{\infty} \frac{\# N_m(C)}{m} z^m\right),$$
see \cite{NogA}. 
This function belongs to the power series ring $\mathbb{Q}[[z]]$. If one prefers, the zeta function can be considered as a function of a complex variable by choosing $z$ to be a sufficiently small complex number.

\begin{them}
Suppose $\map{P}{\widetilde{C}}{C}$ is a $d$-fold ramified covering of finite categories. Then, we have
$$\zeta_{\widetilde{C}}(z)=\zeta_C(z)^d (1-z)^V.$$
\begin{proof}
By Lemma \ref{degenerate}, we have
\begin{eqnarray*}
\zeta_{\widetilde{C}}(z)&=&\exp \bigg(  \sum_{m=1}^{\infty} \frac{\# N_m(\widetilde{C})}{m} z^m\bigg) \\
&=&\exp \bigg( \sum_{x\in \mathrm{Ob}(C)} \bigg( \sum_{\tilde{x}\in P^{-1}(x)} \bigg( \sum_{m=1}^{\infty} \frac{\# (N_m(\widetilde{C})^{\tilde{x}} -\{\mathbf{1}_{\tilde{x}}\})}{m} z^m \\ 
&&+\sum_{m=1}^{\infty}  \frac{1}{m}z^m       \bigg)\bigg)\bigg) \\
&=&\exp  \bigg(\sum_{x\in \mathrm{Ob}(C)} d \bigg( \sum_{m=1}^{\infty} \frac{\# (N_m(C)^x -\{\mathbf{1}_x\})}{m} z^m \bigg) \\ 
&&+\# N_0(\widetilde{C})\sum_{m=1}^{\infty}  \frac{1}{m}z^m       \bigg) \\
&=&\exp\bigg( d\sum_{m=1}^{\infty} \frac{\# N_m(C)}{m} z^m+ (-d\# N_0(C)+\# N_0(\widetilde{C})) \sum_{m=1}^{\infty}  \frac{1}{m}z^m  \bigg).
\end{eqnarray*}
Since $\sum_{m=1}^{\infty}  \frac{1}{m}z^m =-\log(1-z)$, we obtain the result.
\end{proof}
\end{them}

Let us recall the definition of ramified coverings in the sense of Dold \cite{Dol86}.

Let $d$ be a natural number and $X$ be a topological space. We denote by $\mathrm{SP}^d(X)$ the $d$th symmetric power,
$$\mathrm{SP}^d(X)=\prod^d_{i=1} X/\Sigma_d,$$
where $\sum_d$ is the symmetric group of order $d$. Elements $z=[x_1,\dots,x_d]$ of $\mathrm{SP}^d(X)$ are  written as sums, $z=x_1+\cdots+x_d$.

If $\map{\pi}{X}{Y}$ is a continuous map, then a \textit{$d$-inverse} of $\pi$ is a continuous map $\map{t}{Y}{\mathrm{SP}^{d}(X)}$ satisfying the following conditions:
\begin{enumerate}
\item For any $x$ of $X$, $t\circ \pi(x)$ contains $x$.
\item For any $y$ of $Y$, $(\mathrm{SP}^d(\pi))(t(y))=y+\cdots+y$.
\end{enumerate}

When we consider a simplicial set and its geometric realization, we use the same symbols used in \cite{May92}. See Chapter I and III of \cite{May92}.

\begin{them}
Suppose $\map{P}{\widetilde{C}}{C}$ is a $d$-fold ramified covering. Then, the map $\map{BP}{B\widetilde{C}}{BC}$ has a $d$-inverse.
\begin{proof}
Define a map $\map{t}{BC}{\mathrm{SP}^d(B\widetilde{C})}$ as follows: For $$\mathbf{f}=\xymatrix{x_0\ar[r]^{f_1}&x_1\ar[r]^{f_2}&\cdots\ar[r]^{f_{n}}&x_{n}}$$
of $N_n(C)$ such that $\mathbf{f}\not =\mathbf{1}$, $\mathbf{f}$ has $d$ lifts $\tilde{\mathbf{f}}_1,\tilde{\mathbf{f}}_2,\cdots,\tilde{\mathbf{f}}_d$ by Lemma \ref{degenerate}. Define 
$$t|\mathbf{f},u_n|=|\tilde{\mathbf{f}}_1, u_n|+|\tilde{\mathbf{f}}_2, u_n|+\cdots+|\tilde{\mathbf{f}}_d, u_n|$$
for any $u_n$ of $\Delta_n$. If $\mathbf{f}=\mathbf{1}$, there exists unique lift $\tilde{\mathbf{1}}$, and define $$t|\mathbf{1},u_n|=|\tilde{\mathbf{1}}, u_n|+|\tilde{\mathbf{1}}, u_n|+\cdots+|\tilde{\mathbf{1}}, u_n|.$$ We show $t$ is well-defined.

By ignoring the order of lifts, we can exchange to take lifts and face operators, that is, $$\{ \partial_i \tilde{\mathbf{f}}_1,\partial_i \tilde{\mathbf{f}}_2,\dots,\partial_i \tilde{\mathbf{f}}_d \}=\{  \widetilde{(\partial_i \mathbf{f} )}_1,\widetilde{(\partial_i \mathbf{f} )}_2,\dots,\widetilde{(\partial_i \mathbf{f} )}_d \},$$
and this is shown in the following three cases:

If $\mathbf{f}=\mathbf{1}$, then $\partial_i \mathbf{f}=\mathbf{1}$ and $\widetilde{(\partial_i \mathbf{f} )}_j=\mathbf{1}$ for all $j$, so that the two sets consist of $\tilde{\mathbf{f}}$.

If $\mathbf{f}\not =\mathbf{1}$ and $\partial_i \mathbf{f}=\mathbf{1}$, then Lemma \ref{iff} implies $\partial_i \tilde{\mathbf{f}}_j=\tilde{\mathbf{1}}$ for all $j$. Hence, the two sets are equal.

Suppose $\mathbf{f}\not =\mathbf{1}$ and $\partial_i \mathbf{f}\not =\mathbf{1}$. Write each lift $\tilde{\mathbf{f}}_j$ of $\mathbf{f}$ by 
$$\tilde{\mathbf{f}}_j=\xymatrix{\tilde{x}_{j,0}\ar[r]^{\tilde{f}_{j,1}}&\tilde{x}_{j,1}\ar[r]^{\tilde{f}_{j,2}}&\cdots\ar[r]^{\tilde{f}_{j,n}}&\tilde{x}_{j,n}}.$$
Let $k$ be the first number such that $\map{f_k}{x_{k-1}}{x_k}$ of $\mathbf{f}$ is not an identity morphism. If $k<n$, then $x_{j,m} \not =x_{\ell,m}$ for any $j,\ell$ and $m\ge k$. Hence, all $\partial_i \tilde{\mathbf{f}}_j$ are distinct. if $k=n$, then all $\partial_i \tilde{\mathbf{f}}_j$ are also distinct since $i\not = n$. Since $\partial_i \mathbf{f}$ has exactly $d$ lifts, the two sets are equal.

Thus, we have
\begin{eqnarray*}
t|\partial_i \mathbf{f}, u_{n-1}|&=&|\widetilde{(\partial_i \mathbf{f})}_1, u_{n-1}|+|\widetilde{(\partial_i \mathbf{f})}_2, u_{n-1}|+\cdots+|\widetilde{(\partial_i \mathbf{f})}_d, u_{n-1}| \\
&=&|\partial_i \tilde{\mathbf{f}}_1, u_{n-1}|+|\partial_i \tilde{\mathbf{f}}_2, u_{n-1}|+\cdots+|\partial_i \tilde{\mathbf{f}}_d, u_{n-1}| \\
&=&|\tilde{\mathbf{f}}_1, \delta_i u_{n-1}|+| \tilde{\mathbf{f}}_2, \delta_i  u_{n-1}|+\cdots+|\tilde{\mathbf{f}}_d, \delta_i  u_{n-1}| \\
&=&t|\mathbf{f},\delta_i u_{n-1}|.
\end{eqnarray*}
It is clear that $$\{ s_i \tilde{\mathbf{f}}_1,s_i \tilde{\mathbf{f}}_2,\dots,s_i \tilde{\mathbf{f}}_d \}=\{  \widetilde{(s_i \mathbf{f} )}_1,\widetilde{(s_i \mathbf{f} )}_2,\dots,\widetilde{(s_i \mathbf{f} )}_d \}.$$
Hence, we have $$t|s_i \mathbf{f}, u_{n+1}|=t|\mathbf{f}, \sigma_i u_{n+1}|,$$
so $t$ is well-defined.

Next we show that $t$ satisfies the required properties. For any $|\tilde{\mathbf{f}}, u_n|$ of $B\widetilde{C}$, we have
\begin{eqnarray*}
t\circ BP|\tilde{\mathbf{f}}, u_n|&=&t|P(\tilde{\mathbf{f}}), u_n| \\
&=&t|\mathbf{f}, u_n| \\
&=&|\tilde{\mathbf{f}}_1, u_n|+|\tilde{\mathbf{f}}_2, u_n|+\cdots+|\tilde{\mathbf{f}}_d, u_n|,
\end{eqnarray*}
and $\tilde{\mathbf{f}}_i=\tilde{\mathbf{f}}$ for some $i$. For any $|\mathbf{f},u_n|$ of BC, we have
\begin{eqnarray*}
\mathrm{SP}^d(BP)\circ t|\mathbf{f},u_n|&=&\mathrm{SP}^d(BP)( |\tilde{\mathbf{f}}_1, u_n|+|\tilde{\mathbf{f}}_2, u_n|+\cdots+|\tilde{\mathbf{f}}_d, u_n| ) \\
&=&|\mathbf{f},u_n|+|\mathbf{f},u_n|+\cdots+|\mathbf{f},u_n|.
\end{eqnarray*}
The rest of this proof is due to showing continuity of $t$. We define the following map $T$:

$$  \xymatrix{ \displaystyle \coprod_{n\ge 0}N_n(C)\times \Delta_n \ar[dd]^{Q_1} \ar@{.>}[rr]^(0.45)T&&    \displaystyle \prod^d_{i=1}\bigg(\coprod_{n\ge 0}N_n(\widetilde{C})\times \Delta_n \bigg) \ar[d]^{Q_3} \\ &&  \displaystyle\prod^d_{i=1}B\widetilde{C} \ar[d]^{Q_2}\\BC \ar[rr]^t&& \mathrm{SP}^d(B\widetilde{C}),}$$
where each $Q_i$ is the natural projection. If the map $T$ is continuous and makes the diagram commutative, for any open set $U$ of $\mathrm{SP}^d(B\widetilde{C})$, $t^{-1}(U)$ is open in $BC$ if and only if $Q_1^{-1}(t^{-1}(U))$ is open in $\coprod_{n\ge 0}N_n(C)\times \Delta_n $ if and only if $T^{-1}(Q_3^{-1}(Q_2^{-1}(U)))$ is open in $\coprod_{n\ge 0}N_n(C)\times \Delta_n $, so that $t$ is continuous.

Define $T$ as follows: For any $(\mathbf{f},u_n)$ of $N_n(C)\times \Delta_n$, fix an order of its lifts such as $\tilde{\mathbf{f}}_1,\tilde{\mathbf{f}}_2,\dots,\tilde{\mathbf{f}}_d$. Define
$$T(\mathbf{f},u_n)=((\tilde{\mathbf{f}}_1,u_n),(\tilde{\mathbf{f}}_2,u_n),\dots,(\tilde{\mathbf{f}}_d,u_n) ).$$ To show $T$ is continuous, it suffices to show that for any open set $\{\tilde{\mathbf{g}}_i \}\times U_{n_i}$ for $1\le i\le d$, 
\begin{eqnarray}
T^{-1}\bigg( \displaystyle \prod^d_{i=1} \{\tilde{\mathbf{g}}_i \}\times U_{n_i} \bigg) \label{T}
\end{eqnarray}
is open in $\coprod_{n\ge 0} N_n(C)\times \Delta_n$. Note that $N_n(C)$ has a discrete topology. If $P(\tilde{\mathbf{g}}_i)\not =P(\tilde{\mathbf{g}}_j)$ for some $i$ and $j$, the inverse image \eqref{T} is an empty set. Suppose $P(\tilde{\mathbf{g}}_i)=\mathbf{g}$ for all $i$. Then, if $(\tilde{\mathbf{g}}_1,\tilde{\mathbf{g}}_2,\dots,\tilde{\mathbf{g}}_d)$ has a difference order from the one that we fixed, the inverse image \eqref{T} is also an empty set. In the other case, we have
$$T^{-1}\bigg( \prod^d_{i=1} \{\tilde{\mathbf{g}}_i \}\times U_{n_i} \bigg)=\{\mathbf{g}\} \times \bigg( \bigcap^d_{i=1} U_{n_i} \bigg),$$
so that $T$ is continuous. It is easy to show $T$ makes the diagram commutative.

Hence, we conclude $t$ is a $d$-inverse of $BP$.
\end{proof}
\end{them}

\subsection{Examples}\label{Exam}

Let $C$ be a small category. An object $x$ of $C$ is \textit{preinitial} if $T(x)$ consists of only the identity morphism $1_x$.

\begin{defi}
Suppose $C_1$ and $C_2$ are small categories and $x_1$ and $x_2$ are preinitial objects of $C_1$ and $C_2$, respectively. Define the \textit{wedge product} $C_1 \vee C_2$ of $C_1$ and $C_2$ by the following: The set of objects of $C_1\vee C_2$ is
$$\mathrm{Ob}(C_1\vee C_2)=\mathrm{Ob}(C_1)\coprod \mathrm{Ob}(C_2)/\sim,$$
where only $x_1$ and $x_2$ are identified. For objects $x$ and $y$ of $\mathrm{Ob}(C_1\vee C_2)$, the Hom-set is 
$$\mathrm{Hom}_{C_1\vee C_2}(x,y)=  \begin{cases} \mathrm{Hom}_{C_i}(x,y) & \text{if} \ x,y \in \mathrm{Ob}(C_i) \\ \mathrm{Hom}_{C_i} (x_i,y) & \text{if} \ x=x_1\ \text{or}\ x_2, y\in \mathrm{Ob}(C_i) \\ \emptyset &\text{otherwise}. \end{cases}$$
\end{defi}

\begin{exam}
Suppose $C$ is a small category and $x_1$ is a preinitial object. Then, the natural projection from the $d$-tuple wedge product $\bigvee^d_{i=1} C$ of $C$ to $C$ is a $d$-fold ramified covering.
\end{exam}

\begin{exam}
In fact, we can obtain more general result. Suppose $\map{P_i}{\widetilde{C}_i}{C}$ $(1\le i\le d)$ is an unramified covering and $\tilde{x}_i$ is a preinitial object of $\widetilde{C}_i$ such that $P(\tilde{x}_i)=x$ for all $i$. Then, the $d$-tuple wedge product
$$\lmap{\bigvee^d_{i=1} P_i}{\bigvee^d_{i=1} \widetilde{C}_i}{C}$$
is a $d$-fold ramified covering.
\end{exam}

\begin{exam}
Let $\widetilde{C}$ be the following category:
$$
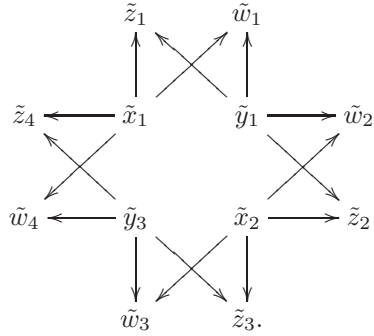
\xymatrix{&\tilde{z}_1&\tilde{w}_1& \\ \tilde{z}_4&\tilde{x}_1\ar[u]\ar[l]\ar[ur]\ar[ld]&\tilde{y}_1 \ar[u]\ar[ul]\ar[r] \ar[rd]&\tilde{w}_2 \\ \tilde{w}_4&\tilde{y}_3 \ar[ul] \ar[l]\ar[d]\ar[rd]&\tilde{x}_2\ar[r] \ar[ur] \ar[d] \ar[dl]&\tilde{z}_2 \\ &\tilde{w}_3&\tilde{z}_3 .&}$$
Let $$C=\xymatrix{x\ar[r]\ar[dr]&z \\y\ar[r]\ar[ur]&w.}$$ Define a functor $\map{P}{\widetilde{C}}{C}$ by eliminating the tildes and the indexes. Then, $P$ is a $4$-fold ramified covering. Their series Euler characteristics are
$$\chi_{\sum}(\widetilde{C})=-4,\ \ \chi_{\sum}(C)=0.$$
Of course, the Riemann-Hurwitz formula holds for $P$ as $V=4$. Their zeta functions are
$$\zeta_{\widetilde{C}}(z)=\frac{1}{(1-z)^{12}}\exp\left(\frac{16z}{1-z}\right),$$
$$\zeta_C(z)=\frac{1}{(1-z)^{4}}\exp\left(\frac{4z}{1-z}\right)$$
by Corollary 2.12 of \cite{NogA}.
\end{exam}


\begin{thebibliography}{AAA99}


\bibitem[BN09]{BN09} M. Baker and S. Norine. Harmonic morphisms and hyperelliptic graphs. \textit{Int. Math. Res. Not. IMRN}, 15: 2914--2955, 2009.

\bibitem [BH99]{BH99} Martin R. Bridson and Andr{\'e} Haefliger. \textit{Metric spaces of non-positive curvature}. Grundlehren der Mathematischen Wissenschaften [Fundamental Principles of Mathematical Sciences], 319, Springer-Verlag, Berlin, 1999.


\bibitem [BL08]{Leib}  C. Berger and T. Leinster. The Euler characteristic of a category as the sum of a divergent series, \textit{Homology, Homotopy Appl., }10(1): 41-51, 2008.

     


\bibitem[CM]{CM} C. Cibils and J. MacQuarrie. Gradings, smash products and Galois coverings of a small category. arXiv:1203.5905v1.

\bibitem[Dol86]{Dol86} A. Dold. Ramified coverings, orbit projections and symmetric powers. \textit{Math. Proc. Cambridge Philos. Soc.}, 99(1): 65--72, 1986.
      
\bibitem[FLS11]{FLS} T. M. Fiore, W. L\"{u}ck and R. Sauer. Finiteness obstructions and Euler characteristics of categories. \textit{Adv. Math}, 226(3), 2371--2469, 2011.



\bibitem [Hat02]{Hat02} A. Hatcher. \textit{Algebraic topology}. Cambridge University Press, Cambridge, 2002.
 

\bibitem [Lei08]{Leia}    T. Leinster. The Euler characteristic of a category. \textit{Doc. Math.}, 13: 21-49, 2008.


\bibitem[MM10]{MM10} B. Malmskog and M. Manes. ``{A}lmost divisibility'' in the {I}hara zeta functions of
              certain ramified covers of {$q+1$}-regular graphs. \textit{Linear Algebra Appl.}, 432(10): 2486--2506, 2010.
     
\bibitem[May92]{May92} J. P. May. \textit{Simplicial objects in algebraic topology.} Chicago Lectures in Mathematics, Reprint of the 1967 original, University of Chicago Press, 1992.

\bibitem [May99]{May99}     J. P. May. \textit{A concise course in algebraic topology} . Chicago Lectures in Mathematics. University of Chicago Press, Chicago, IL, 1999. 


\bibitem [Mor12]{Mor12} M. Morishita. \textit{Knots and primes: An introduction to arithmetic topology}.  Universitext, Springer, London, 2012.



\bibitem [Nog11]{Nog11} K. Noguchi. The Euler characteristic of acyclic categories. \textit{Kyushu Journal of Math.}, 65(1): 85--99, 2011.

\bibitem [Nog]{Nog} K. Noguchi. Euler characteristics of categories and barycentric subdivision. \textit{M\"unster Journal of Math.}, in press. 

\bibitem [NogA]{NogA} K. Noguchi. The zeta function of a finite category. Preprint.

\bibitem[Smi83]{Smi83} L. Smith. Transfer and ramified coverings. \textit{Math. Proc. Cambridge Philos. Soc.}, 93(3): 485--493, 1983.
    
\bibitem[ST96]{ST96} H.M. Stark and A.A. Terras. Zeta functions of finite graphs and coverings. \textit{Adv. in Math.} 121: 124-165, 1996.

\bibitem[Tan]{Tan} K. Tanaka. A model structure on the category of small categories for coverings. \textit{Math. Journal of Okayama Univ.}, in press. 

\bibitem[Ter11]{Ter11} A. A. Terras. \textit{Zeta functions of graphs: A stroll through the garden.} Cambridge Studies in Advanced Mathematics, 128, Cambridge University Press, Cambridge, 2011.
      
\end{thebibliography}
\end{document}